\documentclass[12pt]{article}
\usepackage{amsmath, amsthm, xspace, amssymb}
\usepackage[english]{babel}

\usepackage{graphicx}

\usepackage{ifpdf} 
\ifpdf
\usepackage{epstopdf}
\usepackage{hyperref}
\else
\usepackage[hypertex]{hyperref}
\fi
\usepackage[arrow,matrix, curve]{xy}

\newtheorem{thm}{Theorem}
\newtheorem{lemma}[thm]{Lemma}

\newtheorem{conv}[thm]{Convention}
\theoremstyle{definition}
\newtheorem{rem}[thm]{Remark}
\newtheorem{Def}[thm]{Definition}

\DeclareMathOperator{\id}{id}
\DeclareMathOperator{\End}{End}

\newcommand\erf[1]      {\mbox{(\ref{#1})}}
\newcommand\void[1] {}

\renewcommand\hom  {\mbox{Hom}}
\renewcommand\phi  {\varphi}

\def\be            {\begin{equation}}
\def\bearl         {\begin{array}{l}} 
\def\bearll        {\begin{array}{ll}}
\def\bearlll       {\begin{array}{lll}}

\def\C                     {\ensuremath{\mathcal{C}}\xspace}

\def\CX                                    {\ensuremath{\C^{\mathcal{X}}}\xspace}

\def\congto                {\stackrel{\cong}{\to}}
\def\cover         {\mathcal{F}}
\def\defn         {\stackrel{\mathrm{def}}=}
\def\dimk                        {\mbox{dim}_\Bbbk}
\def\D                          {\mathcal{D}}
\def\ext                   {\ensuremath{{\bf Ext}}\xspace}
\def\ee            {\end{equation}}
\def\eear          {\end{array}}
\def\End                        {\mbox{End}}
\def\Endfun{\mathcal{E}\hspace{-1.6pt}{\it nd}}

\def\Gext                  {\ensuremath{{\bf GExt}}\xspace}

\def\homc                  {\hom_{\C}}

\def\intEnd        {\underline{\End}}
\def\I                     {\mathcal{I}}
\def\la                    {\langle}
\def\M                          {\mathcal{M}}
\def\oti                   {\otimes}

\def\unit          {\mathbf{1}}  

\def\ra                    {\rangle}

\def\tauzz                 {\tau^{\X}}

\def\vectk         {{\mathcal V}\mbox{\sl ect}_\Bbbk}
\def\X                     {\ensuremath{\mathcal{X}}\xspace}
\def\ZZ                    {\ensuremath{\mathbb{Z}/2}\xspace}

\def\B                          {{\bf B}}

\def\bp            {\begin{picture}}
\newcommand\pic[1] {\includegraphics{#1}}
\def\ep            {\end{picture}}

\newcommand\rib[1] {\includegraphics{#1}}

\newcommand{\BIGOP}[1]{\mathop{\mathchoice%
{\raise-0.22em\hbox{\LARGE $#1$}}%
{\raise-0.05em\hbox{\Large $#1$}}{\hbox{\large $#1$}}{#1}}}

\begin{document}

\thispagestyle{empty}
\def\thefootnote{\fnsymbol{footnote}}
\begin{flushright}
   {\sf ZMP-HH/10-15}\\
   {\sf Hamburger$\;$Beitr\"age$\;$zur$\;$Mathematik$\;$Nr.$\;$ 382}\\[2mm]
   June 2010
\end{flushright}
\vskip 2.0em
\begin{center}\Large 
   PERMUTATION MODULAR INVARIANTS FROM MODULAR FUNCTORS

\end{center}\vskip 1.4em
\begin{center}
  Till Barmeier\,$^{\,a}$\footnote{\scriptsize 
  ~Email address: \\
  $~$\hspace*{2.4em}barmeier@math.uni-hamburg.de}
\end{center}

\vskip 3mm

\begin{center}\it$^a$
  Organisationseinheit Mathematik, \ Universit\"at Hamburg\\
  Bereich Algebra und Zahlentheorie\\
  Bundesstra\ss e 55, \ D\,--\,20\,146\, Hamburg
\end{center}

\vskip 2.5em
\begin{abstract} \noindent
For any finite group $G$ with a finite $G$-set $\X$ 
and a modular tensor category $\C$ we construct a part of the
algebraic structure of an associated $G$-equivariant monoidal category: For any
group element $g\in G$ we exhibit the 
module category structure of the $g$-component over the trivial component.
This uses the formalism of permutation equivariant modular functors that was worked out in 
\cite{permMF}.
As an application we show that the corresponding modular invariant partition function
is given by permutation by $g$.
\end{abstract}

\setcounter{footnote}{0} \def\thefootnote{\arabic{footnote}} 

\newpage

\section{Introduction}


The structure of a $G$-equivariant monoidal category has been introduced to understand orbifold models of rational conformal 
field theory with automorphism group $G$ (\cite{kirillov}). A way to describe such categories is given by the 
formalism of $G$-equivariant modular functors (\cite{kpGMF}).
In simple words, a $G$-equivariant modular functor assigns $\vectk$-valued functors to principal $G$-covers.
In \cite{permMF} for
any finite group $G$ a 
$G$-equivariant modular functor $\tau^\X$ was constructed out of a finite $G$-set \X and a modular tensor category $\C$. 
In the present paper we derive certain aspects of the corresponding $G$-equivariant monoidal category \CX. More specifically, 
we exhibit the module categories that are part of \CX.

We will always assume that $\C$ is a $\Bbbk$-linear modular tensor category, where $\Bbbk$ is an algebraically closed field 
of characteristic $0$. In particular, \C has only finitely many simple objects, 
we will write $\I$ for the set of isomorphism classes of simple
objects of $\C$ and $U_i$ with $i\in\I$ for representatives of these. To simplify notation, we agree to drop the tensor 
product symbol for two objects in \C, so we write $AB\equiv A\otimes B$.

We will shortly repeat the construction presented in \cite{permMF}. Let $\tau$ be the \C-extended modular functor that corresponds to the modular category \C (\cite{BK}). 
This is an assignment of a functor
\be
\tau(\Sigma):\C^{\boxtimes A(\Sigma)}\to \vectk
\ee
to any extended surface $\Sigma$, where $A(\Sigma)$ is the set of boundary components of $\Sigma$. After choosing
a total order on $\X$ as an auxiliary datum, the action of
$G$ on $\X$ induces (\cite[Section 2]{permMF}) a functor
\be
\cover_\X:\Gext\to \ext
\ee
from the category of principal $G$-covers of extended surfaces to the category of extended surfaces by taking
the total space of the associated cover. In more detail we put
\be
\cover_\X(P\to M):=\X\times_GP
\ee
for every $G$-cover $(P\to M)$ of an extended surface.
This functor was called the cover functor. Now if $O_g$ is the set of $\la g\ra$-orbits of \X, put
$$
\CX:=\bigoplus_{g\in G}\CX_g
$$
with $\CX_g:=\C^{\boxtimes O_g}$. The assignment
\be
\tauzz(P\to M):=\tau(\cover_\X(P\to M))=\tau(\X\times_GP)
\ee
then gives a $\CX$-extended $G$-equivariant modular functor. For the details of this construction we refer to \cite{permMF}.

In \cite{kpGMF} it was shown that a (genus zero) $\CX$-extended $G$-equivariant modular functor is equivalent to the 
structure of a $G$-equivariant (weakly) fusion category on $\CX$. For the case of $G=\ZZ$ acting on a two-element set
by permutation, the complete set of structure morphisms for this monoidal structure was presented in \cite{permMF} 
by analyzing the
geometric structure of the surfaces $\X\times_{\ZZ}P$ 
using techniques from the Lego-Teichm{\"u}ller-Game (\cite{BKlego}).

For arbitrary finite groups $G$, there is currently no similar result, since the situation is more involved.
The explicit algebraic description of the full monoidal structure of \CX is far out of reach.
As a step towards this, one notices that the axioms of a $G$-equivariant monoidal category imply that the 
summands $\CX_g$ are module categories over the 
monoidal category $\CX_1$. Knowing these structures opens perspectives in two directions: 

On the one hand, a module category $\mathcal{M}$ over a fusion category $\mathcal{D}$ comes with the two $\alpha$-induction functors 
$\alpha^\pm:\mathcal{D}\to\Endfun_{\mathcal{D}}(\mathcal{M})$.
An important quantity is then given by the modular invariant partition function
\be
Z(\mathcal{M}/\mathcal{D})_{i,j}:=\dimk\hom_{\Endfun_{\mathcal{D}}(\mathcal{M})}(\alpha^+_i,\alpha^-_j)\,,
\ee
where $i,j$ label the simple objects of $\mathcal{D}$.

On the other hand, a large part of the structure
of the full $G$-equivariant category is already encoded in the collection of these module categories: In \cite[Section 8]{ENO} it was shown that for a 
fusion category $\mathcal{D}$, a group
homomorphism $c:G\to \mathrm{BrPic}(\mathcal{D})$ from $G$ into the group $\mathrm{BrPic}(\mathcal{D})$ of equivalence classes of invertible module
categories over $\mathcal{D}$ induces two elements in certain cohomology groups of $G$. There exists a structure of a 
$G$-equivariant monoidal category with neutral component $\mathcal{D}$ and as twisted
components representatives of the equivalence classes $c(g)$, if and only if these two obstruction classes vanish. 
Equivalence classes of $G$-equivariant categories based on the homomorphism $c$ 
then form a 
torsor over $H^3(G,\Bbbk^\times)$. Now in our modular functor approach, the existence of a $G$-equivariant monoidal structure
on the system of module categories $\CX_g$ is already ensured by \cite{kpGMF}, hence in our situation both obstruction classes
have to be trivial. Thus the results of the present paper describe the equivalence class of the $G$-equivariant category up
to an element of a torsor over $H^3(G,\Bbbk^\times)$.

In this paper we will exhibit the geometric objects that are relevant for the module category structure and 
derive the corresponding module functors and mixed associativity constraints. Sections \ref{sec:action} and \ref{sec:constraints} culminate in the first
main theorem that describes the structure of a module category over $\CX_1$ on $\CX_g$:
\begin{thm}\label{thm:firstmain}
For any finite group $G$, any finite $G$-set $\X$ and any $g\in G$ the functor
\be
\bearll
\C^{\boxtimes\X}\times \C^{\boxtimes O_g}&\rightarrow \C^{\boxtimes O_g}\\
(A_x)_{x\in\X}\times (M_o)_{o\in O_g}&\mapsto (A_{x_o}A_{g^{-1}x_o}\cdots A_{g^{-|o|+1}x_o}M_o)_{o\in O_g}
\eear
\ee
with $x_o$ the smallest element in the $\la g\ra$-orbit $o$, together with the associativity constraints
\be
\Psi_{A,B,M}=(\psi^o_{(A_x)_{x\in o},(B_x)_{x\in o},M})_{o\in O_g}\,,
\ee
where the morphisms $\psi^o$ contain only braiding morphisms as in equation \erf{eqn:psiexample}, endows $\C^{\boxtimes O_g}$ with the structure of a module category over the tensor category $\C^{\boxtimes \X}$.
\end{thm}
In section \ref{sec:invariants} we prove the following theorem about the modular invariant 
matrix $Z(\CX_g/\CX_1)$ of the module category $\CX_g$ over $\CX_1$:

\begin{thm}
The modular matrix $Z(\CX_g/\CX_1)$ for the module category described in theorem \ref{thm:firstmain} reads
\be
Z(\CX_g/\CX_1)_{\bar \imath,\bar \jmath}=\delta_{\bar \jmath,g\bar\imath}
\ee
where $\bar\imath, \bar\jmath\in\I^\X$ label the simple objects of $\CX_1=\C^{\boxtimes\X}$ and 
$g\bar\imath$ is the multi-index $\bar\imath$ permuted by the action of the group element $g\in G$.
\end{thm}

\subsubsection*{Acknowledgements}
The author would like to thank Thomas Nikolaus and Christoph Schweigert for helpful discussions
and comments on the draft.
This work was supported by the DFG Priority Program SPP 1388
``Representation theory''.

\section{Module categories from $G$-equivariant modular functors}\label{sec:main}
\subsection{Preliminaries}
Recall from \cite{permMF} that there is the structure of a $G$-equivariant monoidal category on 
$\bigoplus_{g\in G}\CX_g$ with $\CX_1=\C^{\boxtimes\X}$ and $\CX_g=\C^{\boxtimes O_g}$. We first want to find
the monoidal structure on the neutral component $\C^{\boxtimes\X}$ and then for every $g\in G$ the module action functor
$$
\C^{\boxtimes\X}\times\C^{\boxtimes O_g}\to\C^{\boxtimes O_g}\,\,.
$$
We first briefly turn our attention to the monoidal structure on $\CX_1$. In this case, all relevant $G$-covers of extended surfaces are trivial covers. Since the cover functor $\cover_\X$ maps trivial covers to disjoint unions of copies of the base space, the monoidal structure on $\CX_1$ is found by evaluating the modular functor $\tau$ on disjoint unions of standard $n$-pointed spheres for appropriate $n$. The occurring marking graphs are in all cases the standard marking graphs on $S_n$. Now the following lemma is an easy observation:
\begin{lemma}
The weakly ribbon structure on $\CX_1=\C^{\boxtimes \X}$ induced by the $G$-equivariant modular functor $\tauzz$ is
ribbon and is equivalent to the standard ribbon structure on $\C^{\boxtimes \X}$. The tensoriality constraints of the permutation action of $G$ on $\C^{\boxtimes\X}$ are identities.
\end{lemma}

\subsection{The module action functor}\label{sec:action}

Recall from \cite[Section 3]{prince} or \cite[Section 2.2.2]{permMF} the definition of the standard block $S_n(g_1,\dots,g_n;h_1,\dots,h_n)$ as 
explicit principal $G$-bundles with $n$ marked points over the standard sphere $S_n$. The standard sphere was itself introduced 
in \cite[Section 2.3]{BKlego} as the Riemann sphere $\overline{\mathbb{C}}$ with $n$ holes removed around the natural numbers
$1,\dots,n$.

The action of $\CX_1$ on $\CX_g$ is found by evaluating the $G$-equivariant modular functor $\tau^\X$ on the principal 
$G$-cover $(S_3(g^{-1},1,g;1,1,1)\to S_3)$ of the three-holed sphere. 
In \cite[Lemma 11]{permMF} the connected components of the total space of the associated bundle
$E_{g_1;g_2}:=\cover_\X(S_3(g_1,g_2,(g_1g_2)^{-1};1,1,1)\to S_3)$ where fully described:

\begin{lemma}\label{lem:poporbits}\mbox{}\\[-2em]
\begin{enumerate}
\item[(i)]
There is a natural bijection between the connected 
components of 
$$\cover_\X(S_3(g_1,g_2,(g_1g_2)^{-1};1,1,1)\to S_3)=E_{g_1;g_2}$$ 
and orbits of the $G$-set $\X$ under the action of the subgroup 
$\langle g_1,g_2\rangle\subset G$ of $G$ generated by 
the elements $g_1$ and $g_2$.
\item[(ii)]
The restriction of $E_{g_1;g_2}$ to the boundary with monodromy $g_1$ is diffeomorphic to 
$E_{g_1^{-1}}:=\mathbb{R}\times\X/(t+2\pi,x)\sim (t,g_1x)$
and similarly for the other boundaries.
Let $o$ be a $\langle g_1,g_2\rangle$-orbit of $\X$ and write $E_{g_1;g_2}^o$ for the connected component of
$E_{g_1;g_2}$ corresponding to the orbit $o$. The boundary components of $E_{g_1;g_2}^o$ correspond to precisely those orbits of the
cyclic subgroups $\langle g_1\rangle$, $\langle g_2\rangle$ and $\langle g_1g_2\rangle$ that are contained in the orbit $o$
of the group $\langle g_1,g_2\rangle$.

\item[(iii)] In particular, the number of sheets
of the cover $E_{g_1;g_2}^o\to S_3$ is $|o|$.
\end{enumerate}
\end{lemma}

By \cite[Lemma 12]{permMF} the genus of the relevant surface $E_{g^{-1};1}$ is zero.

From now on we restrict our attention to the connected components of $E_{g^{-1};1}$ which by lemma \ref{lem:poporbits}
is the same as fixing a $\la g\ra$-orbit $o$.
When we view the corresponding component $E_{g^{-1};1}^o$ as the total space of an $|o|$-fold cover of $S_3$, it has one 
boundary component over the first and third boundary of $S_3$ respectively and $|o|$ boundary components over the second 
boundary of $S_3$. 

In the definition of the module action functor
$$
\C^{\boxtimes \X}\times \C^{\boxtimes O_g}\to \C^{\boxtimes O_g}
$$
the connected component $E_{g^{-1};1}^o$ will give a contribution
$$
\C^{\boxtimes o}\times\C\to\C\,\,;
$$
all these contributions are then factor-wise combined to give the full functor.
By describing this contribution for all $\la g\ra$-orbits separately, we get the full module action functor.

In this way we can restrict ourselves to one connected component. Hence we adopt:
\begin{conv}\label{conv:cyclic}
$G$ is a cyclic group with generator $g$. The ordered set $\X$ has a single $G$-orbit, its smallest element is $x_0\in\X$. Let $n:=|\X|$.
\end{conv}
This convention allows us to simplify notation. In particular we have $\X=\{g^kx_0| k=0,\dots,n-1\}$. 
Note that there is no further assumption on the action of $G$ on $\X$. So $\X=\{x_0\}$ a one-element-set and hence $n=1$
is possible. So from now on we only consider a functor $\CX\times\C\to\C$.

We will write objects in $\CX_1$ as $(A_x)_{x\in\X}$ and sometimes use the abbreviation $(A_x)$. The order of \X induces an order 
on the factors in $\CX_1=\C^{\boxtimes\X}$. When we draw pictures, we will occasionally write $A^{k}$ for the factor $A_{g^kx_0}$ to provide a clearer view of the drawing.

Now let $(A_x)_{x\in\X}$ be an object of $\C^{\boxtimes \X}$ and $M$ an object in $\C$. The tensor product $(A_x)_{x\in \X}\otimes M$ is defined to be the object of $\C$ that represents the functor
\be\label{eqn:repfunctor}
\bearll
\C&\to \vectk\\
T&\mapsto \tauzz(S_3(g^{-1},1,g;1,1,1)\to S_3;T,(A_x)_{x\in\X},M)=\tau(E_{g^{-1};1};T,A_x,M)\,\,.
\eear
\ee
Since the oriented manifold $E_{g^{-1};1}$ is of genus zero and has $n+2$ boundary components, there is a diffeomorphism $E_{g^{-1};1}\cong S_{n+2}$. The choice of such a diffeomorphism induces a natural isomorphism
\be
\tau(E_{g^{-1};1};T,A_x,M)\congto \tau(S_{n+2};T,A_x,M)\defn\homc(\unit,T\otimes(\bigotimes_{x\in\X}A_x)\otimes M)\,\,,
\ee
which then gives a choice $(\bigotimes_{x\in\X}A_x)\otimes M$ of the object representing the functor \erf{eqn:repfunctor}. 

To find an appropriate diffeomorphism, we will draw a marking graph on the surface $E_{g^{-1};1}$. This is most 
conveniently done by viewing $E_{g^{-1};1}$ as the total space of a cover over $S_3$ and then lifting paths 
in $S_3$ to $E_{g^{-1};1}$.

The marking will have one vertex for every boundary component of $E_{g^{-1};1}$
and a single vertex that is connected by an edge to every boundary component. We call this vertex the internal vertex.

Before we turn to concrete graphs, we first have a look at the lifting properties of $E_{g^{-1};1}$. 
Consider the following path in $S_3$ that turns clockwise around the third boundary circle and 
has its starting point $p$ in the lower half-plane:

\be
\raisebox{20pt}{
\bp(146,146)
\put(0,0){\pic{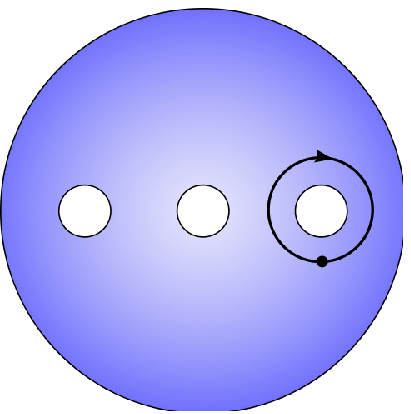}}
\put(22,55){\scriptsize $T$}
\put(55,55){\scriptsize $A$}
\put(89,55){\scriptsize $M$}
\ep
}
\ee

Now we lift this path to $E_{g^{-1};1}$ with starting point $[x,p]$ and find that its end point is $[g^{-1}x,p]$. We will later
use this type of path in order to connect the boundary components over the second boundary circle of $S_3$ 
to our marking graph.

As a first step to finding a marking graph on the cover, we find a path that connects the first and third boundary component of $E_{g^{-1};1}$. 
By definition of the standard block $S_3(g^{-1},1,g;1,1,1)$, the marked point on its first boundary component is 
$p_1=[1-\frac{i}{3},1_G]$ and similarly for the third boundary component $p_3=[3-\frac{i}{3},1_G]$. 
By \cite[Secion 2.2]{permMF} this gives the marked points $[x_0,p_1], [x_0,p_3]\in E_{g^{-1};1}$. 
Now consider the following path in $S_3$:

\be\label{pic:first-third-circle}
\raisebox{20pt}{
\bp(146,146)
\put(0,0){\pic{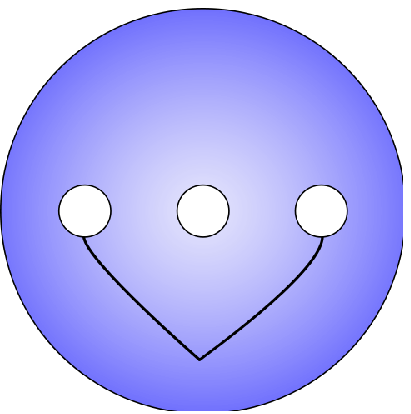}}
\put(22,55){\scriptsize $T$}
\put(55,55){\scriptsize $A$}
\put(89,55){\scriptsize $M$}
\ep
}
\ee

Lifting this path to $E_{g^{-1};1}$ such that the vertex at the first boundary circle is lifted to $[x_0,p_1]$ gives a path that connects $[x_0,p_1]$ and  $[x_0,p_3]$ in $E_{g^{-1};1}$. The lift of the point where the path has the sharp bend
is denoted by $\hat p$ and will serve as internal vertex of the marking graph.

Over the second boundary circle of $S_3$ the surface $E_{g^{-1};1}$ has one boundary component in every sheet with marked points $[x,p_2]$ for every $x\in\X$. By convention \ref{conv:cyclic} we have $x=g^lx_0$ for some $l$. Since $G$ and $\X$ are
finite, we can choose that $x=g^{-k}x_0$ for some $k=0,\dots,n-1$. Now consider the following path in $S_3$ that winds $k$ times clockwise around the third boundary circle:

\be
\raisebox{20pt}{
\bp(146,146)
\put(0,0){\pic{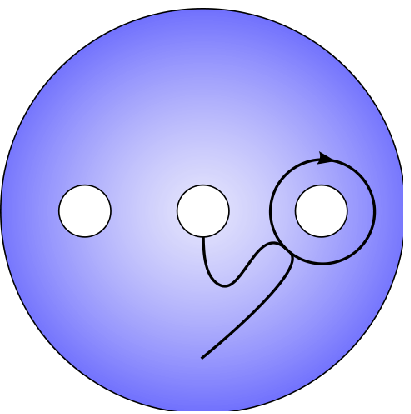}}
\put(22,55){\scriptsize $T$}
\put(55,55){\scriptsize $A$}
\put(89,55){\scriptsize $M$}
\ep
}
\ee
This should be read as follows: The path starts in the lower half-plane, moves near the third boundary component, winds
$k$ times clockwise around it and then connects to the marked point of the second boundary component.

As radius of the circular part of this path we choose $r_k=\frac13+\frac{1}{10}(1-\frac{k-1}{n-1})$.
The lift of this path to $E_{g^{-1};1}$ with starting point $\hat p$ then has end point $[g^{-k}x_0,p_2]$. 
Now we draw the lift of this path on $E_{g^{-1};1}$ for every $k=1\dots n-1$. The assumption on the radius $r_k$
ensures that the graph has no self-intersections, since the radius decreases with increasing $k$. Any other choice of radius
with this property gives a homotopic path. 
For $k=0$ we connect the internal vertex $\hat p$ and $[x_0,p_2]$ with a straight line. This is equivalent to the paths constructed above,
as the straight line is homotopic to the path that winds zero times around the third boundary circle.
In the cover $E_{g^{-1};1}$, this path does not intersect with the lift of the path \erf{pic:first-third-circle}, that connects
the internal vertex to the boundary components labeled by $T$ and $M$.

We finally obtain a marking on $E_{g^{-1};1}$ that connects all marked points on all boundary components. Now we get
a diffeomorphism to $S_{n+2}$ by moving the boundary components of $E_{g^{-1};1}$ along the marking graph and obtain

\be
\raisebox{20pt}{
\bp(146,146)
\put(0,0){\pic{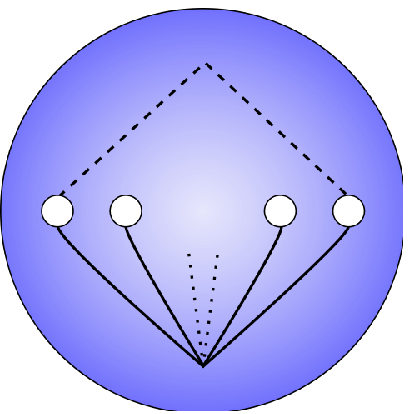}}
\put(4,55){\scriptsize $T$}
\put(54,55){\scriptsize $\dots$}
\put(105,55){\scriptsize $M$}
\put(32,65){\scriptsize $A^{0}$}
\put(67,64){\scriptsize $A^{-n+1}$}
\ep
}
\ee
where the dashed line marks multiple self-intersections of the immersion of the surface into three-dimensional space.
As a non-embedded manifold, this is diffeomorphic to

\be
\raisebox{20pt}{
\bp(146,146)
\put(0,0){\pic{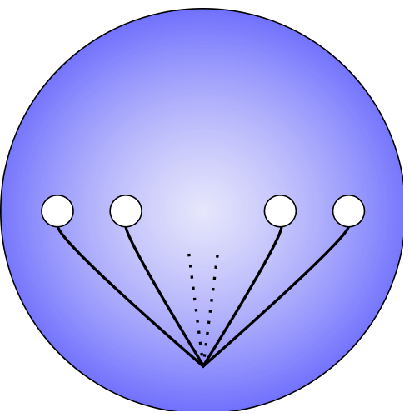}}
\put(15,65){\scriptsize $T$}
\put(54,55){\scriptsize $\dots$}
\put(97,65){\scriptsize $M$}
\put(32,65){\scriptsize $A^{0}$}
\put(73,65){\scriptsize $A^{-n+1}$}
\ep
}
\ee

This diffeomorphism induces an isomorphism
\be
\bearll
\tau(E_{g^{-1};1};T,(A_x)_{x\in\X},M)&\congto \tau(S_{n+2};T,A_{x_0}, A_{g^{-1}x_0},\dots,A_{g^{-n+1}x_0},M)\\
&\defn\homc(\unit,TA_{x_0}\dots A_{g^{-n+1}x_0}M)\,.
\eear
\ee
This shows
\begin{lemma}
The functor 
\be
\bearll
\C&\to \vectk\\
T&\mapsto \tauzz(S_3(g^{-1},1,g;1,1,1)\to S_3;T,(A_x)_{x\in\X},M)\,.
\eear
\ee
is represented by the object
\be
(A_x)_{x\in\X}\otimes M:=A_{x_0}A_{g^{-1}x_0}\dots A_{g^{-n+1}x_0}M\,,
\ee
which serves as a module action functor
\be
\C^{\boxtimes \X}\times\C\to\C\,.
\ee
\end{lemma}

A comment on the order of the objects $A_x$ is due: In $\CX_1=\C^{\boxtimes\X}$ the factors are ordered by the order of
the $G$-set $\X$. In the the module action, the order is by decreasing powers of the generator $g$.

\subsection{Associativity constraints}\label{sec:constraints}

We now turn to the question of finding associativity constraints 
$$
\psi_{(A_x),(B_x),M}:((A_x)\otimes (B_x))\otimes M\to (A_x)\otimes ((B_x)\otimes M)
$$
The general procedure of reading off associativity constraints from $G$-equivariant modular functors was described in
\cite{permMF} in greater detail. 
When dealing with arbitrary groups, even with the restriction to cyclic groups in convention \ref{conv:cyclic}, the 
analysis of covers of the $4$-punctured sphere is rather complicated. The main downside is that we are no longer
able to draw marking graphs on the manifolds $\X\times_GP$ themselves, but have to view them as total spaces of covers
over $S_4$
and lift paths, as in the definition of the module action. This section is very technical, its results are summarized in
theorem \ref{thm:firstmain2}.

We again restrict ourselves to convention \ref{conv:cyclic} and proceed as in \cite[Section 4]{permMF}:
\begin{enumerate}
\item Determine the two marking graphs on 
$\X\times_GS_4(g^{-1},1,1,g;1,1,1,1)$ induced on the cover by 
cutting $S_4$ in the two ways determined by associativity and by our definition of the module action. 
Denote the marking graph representing $((A_x)\otimes (B_x))\otimes M$ by $m_1$ and the marking graph representing $(A_x)\otimes ((B_x)\otimes M)$ by $m_2$.
\item Transform the surface $\X\times_GS_4(g^{-1},1,1,g;1,1,1,1)$ 
with the marking graph $m_1$ to the standard sphere $S_{2n+2}$
in the way prescribed by the marking $m_2$.
\item This yields a marking graph on $S_{2n+2}$. Determine the Lego Teichm{\"u}ller Game (LTG) moves that transform this graph into the standard marking graph on $S_{2n+2}$ and translate these LTG-moves into morphisms in \C.
\end{enumerate}

We turn to the first cutting procedure, that represents the tensor product
$((A_x)\otimes (B_x))\otimes M$. In this case the surface $S_4(g^{-1},1,1,g;1,1,1,1)$ is cut into
a trivial 
$G$-cover $S_3(1,1,1;1,1,1)$ representing the tensor product
$(A_x)\otimes (B_x)$ and the $G$-cover $S_3(g^{-1},1,g;1,1,1)$ representing the product $(C_x)\otimes M$ 
with $(C_x)=(A_x)\otimes (B_x)=(A_xB_x)$. As the cover functor respects gluing, we can analyze the process by 
considering $\X\times_GS_4(g^{-1},1,1,g;1,1,1,1)$ and the respective associated covers over $S_3$. 
We analyze the resulting marking graph on $\X\times_GS_4(g^{-1},1,1,g;1,1,1,1)$ by considering paths in the base $S_4$
and lifting these to the covers. 
For the first and fourth boundary component we get a lift of the path

\be
\raisebox{20pt}{
\bp(146,146)
\put(0,0){\pic{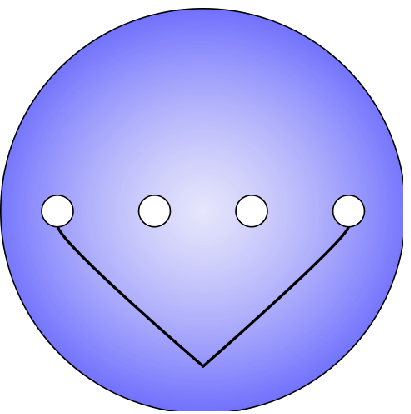}}
\put(14,65){\scriptsize $T$}
\put(97,65){\scriptsize $M$}
\put(41,65){\scriptsize $A$}
\put(69,65){\scriptsize $B$}
\ep
}
\ee
which is again lifted into the sheet corresponding the the generator $x_0\in\X$. Now for the paths that connect to the
boundary components over the second and third boundary circle we consider every sheet separately. For the boundary components in the sheet that corresponds to $x=g^{-k}x_0\in \X$, gluing gives an edge which is a lift of

\be
\raisebox{20pt}{
\bp(300,146)
\put(0,0){\pic{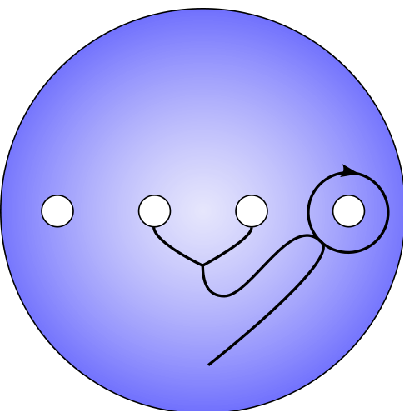}}
\put(14,65){\scriptsize $T$}
\put(97,75){\scriptsize $M$}
\put(41,65){\scriptsize $A$}
\put(69,65){\scriptsize $B$}
\put(145,50){$\stackrel{\cdot}{=}$}
\put(180,0){\pic{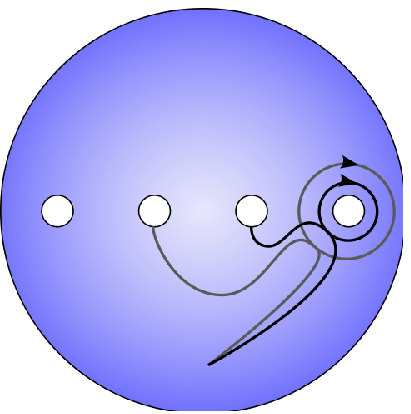}}
\put(180,0){%
\put(14,65){\scriptsize $T$}
\put(97,75){\scriptsize $M$}
\put(41,65){\scriptsize $A$}
\put(69,65){\scriptsize $B$}
}
\ep
}
\ee
where the path turns $k$ times around the fourth boundary circle. In the second picture we contracted the graph along the
factorizing link and draw the edges that connect to different boundary circles in different colors to avoid confusion. 
Note that the intersections of this path do not give intersections in the total space of the cover as turning around the fourth boundary 
circle lifts to a path in the total space that connects different sheets.

In the second cutting procedure the $G$-cover $S_4(g^{-1},1,1,g;1,1,1,1)$ is cut into two $G$-covers 
$S_3(g^{-1},1,g;1,1,1)$, one representing the tensor product
$(B_x)\otimes M$ and one representing the tensor product $(A_x)\otimes N$ with $N=(B_x)\otimes M$. In this case for the 
first and fourth boundary component we again get the path

\be
\raisebox{20pt}{
\bp(146,146)
\put(0,0){\pic{ass-first-fourth}}
\put(14,65){\scriptsize $T$}
\put(97,65){\scriptsize $M$}
\put(41,65){\scriptsize $A$}
\put(69,65){\scriptsize $B$}
\ep
}
\ee
as in the first cutting procedure. For the second and third boundary circle, we again consider all sheets separately; for the sheet corresponding to $g^{-k}x_0\in\X$ we get

\be
\raisebox{20pt}{
\bp(146,146)
\put(0,0){\pic{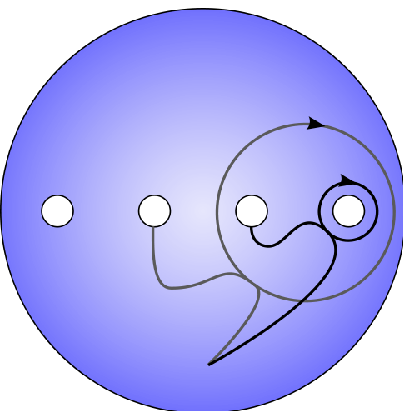}}
\put(14,65){\scriptsize $T$}
\put(97,69){\scriptsize $M$}
\put(41,65){\scriptsize $A$}
\put(69,65){\scriptsize $B$}
\ep
}
\ee
where the path connecting to the boundary component labeled by $B$ turns $k$ times around the boundary component labeled by
$M$, whereas the path connecting to the $A$-boundary turns $k$ times around both the $B$- and the $M$-boundary. Again we use different colors to distinguish the edges.

Now we use the diffeomorphism given by the second marking to transform the manifold into the standard sphere $S_{n+2}$. 
To find the image of the first marking on $S_{2n+2}$, we pick two boundary components of
$\X\times_GS_4(g^{-1},1,1,g;1,1,1,1)$ and see how the corresponding edges of the marking behave relative to each other
while applying the diffeomorphism to $S_{2n+2}$.

The first thing to notice is that the edges connecting to the first and last boundary components do not interfere with
the edges of any other boundary in the application of the diffeomorphism. On $S_{2n+2}$ this just gives

\be
\raisebox{20pt}{
\bp(146,146)
\put(0,0){\pic{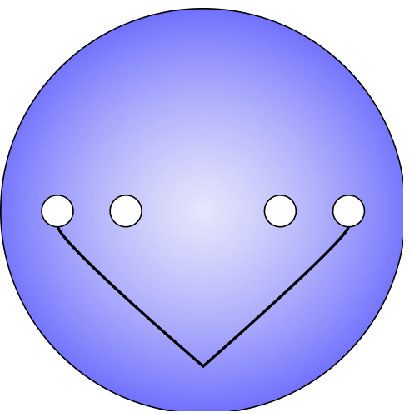}}
\put(14,65){\scriptsize $T$}
\put(97,65){\scriptsize $M$}
\put(53,57){\scriptsize $\dots$}
\ep
}
\ee

We now turn to the boundary components over the second and third boundary circle of $S_4$. When applying the diffeomorphism to $S_{2n+2}$,
all boundary components of $\X\times_GS_4(g^{-1},1,1,g;1,1,1,1)$ are moved simultaneously. When checking the relative 
behavior of two boundary components we will freely move these boundary components and the corresponding edges of the marking.
If the edge comes near any other boundary component over the second and third boundary circle of $S_4$, we will
assume that this component is already moved out of the way or is moved at the same time. This allows us to move the edge
over the other boundary components over the second and third boundary circle of $S_4$. So the following analysis can be 
seen as a kind of recursive algorithm to transform $\X\times_GS_4(g^{-1},1,1,g;1,1,1,1)$ into $S_{2n+2}$. The reader should
always be aware of this procedure and should check that the simultaneous movement indeed justifies this process. 

We will now distinct all possible choices of boundary components over the second and third boundary circle of $S_4$.

\begin{itemize}
\item We start by comparing two boundary components over the third boundary circle of $S_4$, i.e. two boundary components labeled
by $B_{g^{-k}x_0}$ and $B_{g^{-l}x_0}$, where without loss of generality we assume $l> k$. In the markings obtained from the two gluing procedures, the
edges corresponding to $B_{g^{-k}x_0}$ and $B_{g^{-l}x_0}$ are lifts of
\be
\raisebox{20pt}{
\bp(300,146)
\put(0,0){\pic{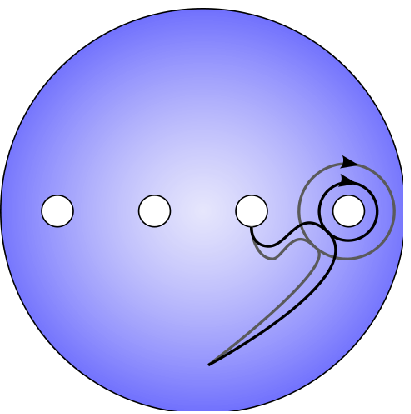}}
\put(14,65){\scriptsize $T$}
\put(97,75){\scriptsize $M$}
\put(41,65){\scriptsize $A$}
\put(69,65){\scriptsize $B$}
\put(180,0){\pic{assltg-third-third}}
\put(180,0){
\put(14,65){\scriptsize $T$}
\put(97,75){\scriptsize $M$}
\put(41,65){\scriptsize $A$}
\put(69,65){\scriptsize $B$}
}
\ep
}
\ee
where the first picture shows the gluing for $((A_x)\otimes(B_x))\otimes M$ and the second picture for 
$(A_x)\otimes ((B_x)\otimes M)$ as explained above. The darker line connects to the boundary component of $B_{g^{-l}x_0}$ and turns $l$ times around the
fourth boundary circle
while the lighter line performs $k$ turns and connects to the $B_{g^{-k}x_0}$-boundary. Obviously both markings coincide, hence on $S_{2n+2}$ we get

\be
\raisebox{20pt}{
\bp(146,146)
\put(0,0){\pic{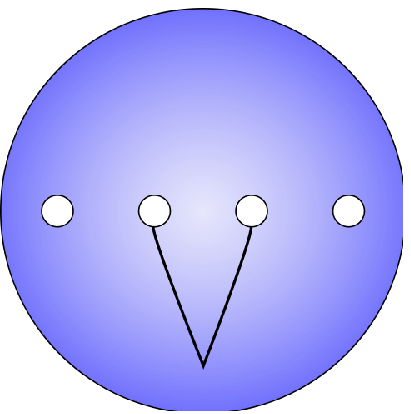}}
\put(42,65){\scriptsize $B^{-k}$}
\put(69,65){\scriptsize $B^{-l}$}
\put(54,57){\scriptsize $\dots$}
\put(25,57){\scriptsize $\dots$}
\put(81,57){\scriptsize $\dots$}
\ep
}
\ee

\item For two boundaries over the second circle of $S_4$ labeled by $A_{g^{-k}x_0}$ and $A_{g^{-l}x_0}$ with $l>k$ we get a similar 
picture. In the gluing procedures described above we obtain edges for the respective boundaries that are lifts of

\be
\raisebox{20pt}{
\bp(300,146)
\put(0,0){\pic{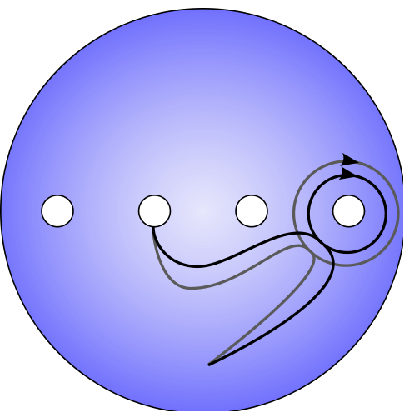}}
\put(14,65){\scriptsize $T$}
\put(97,75){\scriptsize $M$}
\put(41,65){\scriptsize $A$}
\put(69,65){\scriptsize $B$}
\put(180,0){\pic{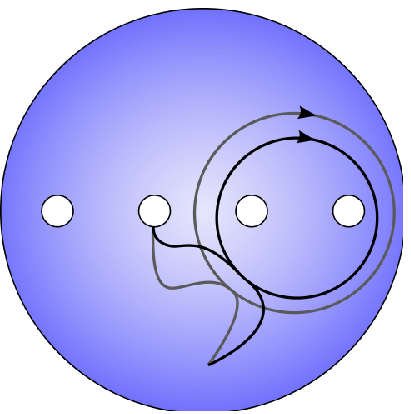}}
\put(180,0){
\put(14,65){\scriptsize $T$}
\put(93,65){\scriptsize $M$}
\put(41,65){\scriptsize $A$}
\put(70,65){\scriptsize $B$}
}
\ep
}
\ee

As we assume that the boundary components over the third circle are already moved out of the way, both edges can be transformed into each other,
hence on $S_{2n+2}$ we get:

\be
\raisebox{20pt}{
\bp(146,146)
\put(0,0){\pic{s2n-third-third}}
\put(41,65){\scriptsize $A^{-k}$}
\put(70,65){\scriptsize $A^{-l}$}
\put(54,57){\scriptsize $\dots$}
\put(25,57){\scriptsize $\dots$}
\put(81,57){\scriptsize $\dots$}
\ep
}
\ee

\item Now we turn to the more complicated situations. For $k\le l$ we compare the boundary components labeled 
by $A_{g^{-k}x_0}$ and $B_{g^{-l}x_0}$. The gluing procedures give us markings where the relevant edges are lifts of 

\be
\raisebox{20pt}{
\bp(300,146)
\put(0,0){\pic{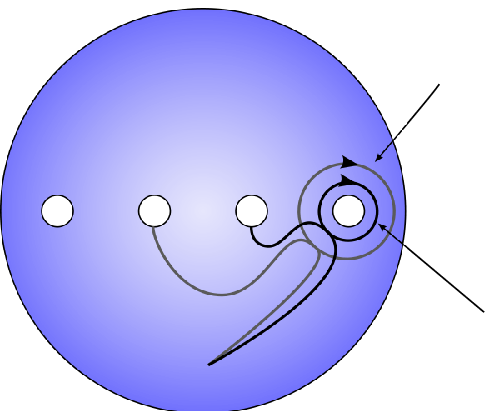}}
\put(14,65){\scriptsize $T$}
\put(97,75){\scriptsize $M$}
\put(41,65){\scriptsize $A$}
\put(69,65){\scriptsize $B$}
\put(112,100){\scriptsize $k$ turns}
\put(129,20){\scriptsize $l$ turns}
\put(180,0){\pic{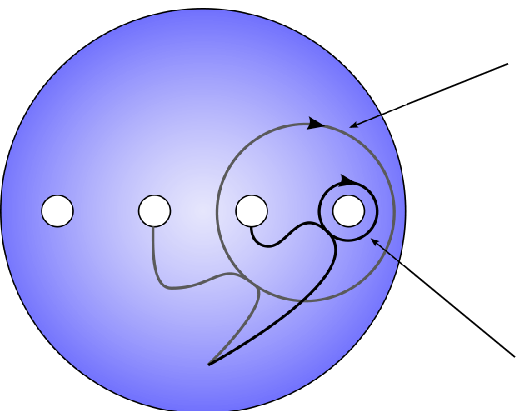}}
\put(180,0){
\put(14,65){\scriptsize $T$}
\put(97,69){\scriptsize $M$}
\put(41,65){\scriptsize $A$}
\put(69,65){\scriptsize $B$}
\put(132,105){\scriptsize $k$ turns}
\put(133,10){\scriptsize $l$ turns}
}
\ep
}
\ee
We see that under the assumption that other boundaries are already moved out of the way, again both edges coincide for $k\le l$. Hence when 
applying the diffeomorphism to $S_{2n+2}$ given by the second marking, we get

\be
\raisebox{20pt}{
\bp(146,146)
\put(0,0){\pic{s2n-third-third}}
\put(39,65){\scriptsize $A^{-k}$}
\put(69,65){\scriptsize $B^{-l}$}
\put(54,57){\scriptsize $\dots$}
\put(25,57){\scriptsize $\dots$}
\put(81,57){\scriptsize $\dots$}
\ep
}
\ee
on $S_{2n+2}$.

\item Finally we compare the edges corresponding to boundary components labeled by $A_{g^{-k}x_0}$ and $B_{g^{-l}x_0}$ with $k>l$. 
In this case the edges in the markings are lifts of

\be
\raisebox{20pt}{
\bp(300,146)
\put(0,0){\pic{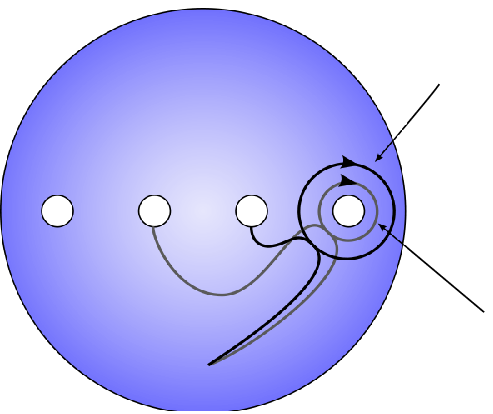}}
\put(14,65){\scriptsize $T$}
\put(97,75){\scriptsize $M$}
\put(41,65){\scriptsize $A$}
\put(69,65){\scriptsize $B$}
\put(112,100){\scriptsize $l$ turns}
\put(129,20){\scriptsize $k$ turns}
\put(180,0){\pic{ass2-second-third-arrows}}
\put(180,0){
\put(14,65){\scriptsize $T$}
\put(97,69){\scriptsize $M$}
\put(41,65){\scriptsize $A$}
\put(69,65){\scriptsize $B$}
\put(132,105){\scriptsize $k$ turns}
\put(133,10){\scriptsize $l$ turns}
}
\ep
}
\ee
Observe that in the first picture the path connecting the internal vertex to the $A_{g^{-k}x_0}$-boundary turns around the fourth circle of $S_4$
with a smaller radius than the path connecting to the $B_{g^{-l}x_0}$-boundary, since $k>l$. Also check that all crossings in
the paths in $S_4$ do not give self-intersections in $\X\times_GS_4(g^{-1},1,1,g;1,1,1,1)$ as the crossing sections of
the paths lift to different sheets. Now we carefully apply the diffeomorphism represented by the second marking. It
instructs us to turn the $A_{g^{-k}x_0}$-boundary $k$ times around the third and fourth boundary circle and the 
$B_{g^{-l}x_0}$-boundary $l$ times around the fourth boundary circle. As a first step we turn the $A_{g^{-k}x_0}$-boundary $(k-l)$
times around. This transforms the first marking into

\be
\raisebox{20pt}{
\bp(146,146)
\put(0,0){\pic{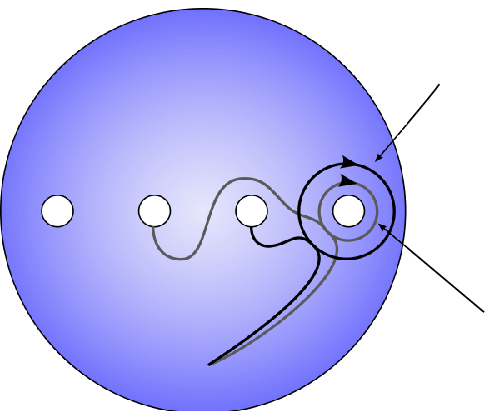}}
\put(14,65){\scriptsize $T$}
\put(97,75){\scriptsize $M$}
\put(41,65){\scriptsize $A$}
\put(69,71){\scriptsize $B$}
\put(112,100){\scriptsize $l$ turns}
\put(129,20){\scriptsize $l$ turns}
\ep
}
\ee
Here both paths wind $l$ times around the fourth boundary circle. When we turned around the $A_{g^{-k}x_0}$-boundary, the edge
connecting to it always passed along the $B_{g^{-l}x_0}$-boundary as they were lying in different sheets. Now we turn both
boundaries around the fourth circle $l$ times simultaneously and finally end up with the marking

\be
\raisebox{20pt}{
\bp(146,146)
\put(0,-2){\pic{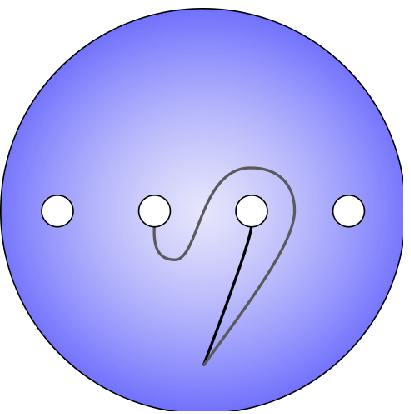}}
\put(42,65){\scriptsize $A^{-k}$}
\put(69,73){\scriptsize $B^{-l}$}
\put(54,55){\scriptsize $\dots$}
\put(25,55){\scriptsize $\dots$}
\put(81,55){\scriptsize $\dots$}
\ep
}
\ee
on $S_{2n+2}$.
\end{itemize}
This describes the relative position of all pairs of edges of the marking we obtain on $S_{2n+2}$. An example of the final
marking in the case of $n=4$ is depicted in 

\be
\raisebox{20pt}{
\bp(280,146)
\put(0,0){\pic{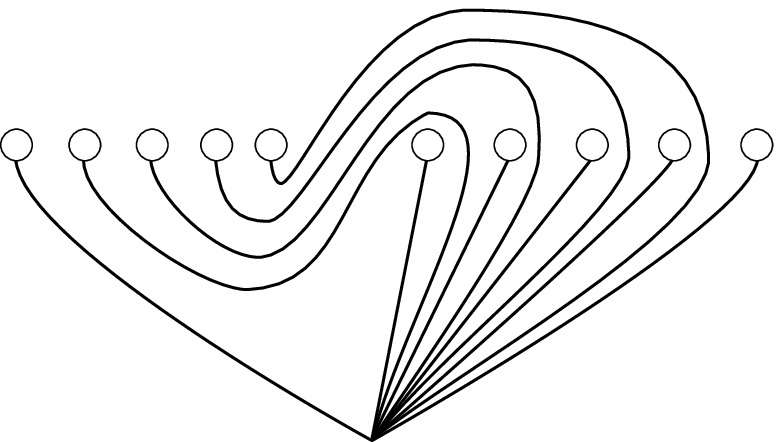}}
\put(2,95){\scriptsize $T$}
\put(215,95){\scriptsize $M$}
\put(39,95){\scriptsize $A's$}
\put(162,95){\scriptsize $B's$}
\ep
}
\ee

In the general case, the final marking on
$S_{2n+2}$ now has straight lines that connect the internal vertex to the $T$-, the $M$-, the $B_x$- and to the $A_{x_0}$-boundaries. 
The edge that connects the internal vertex to the $A_{g^{-k}x_0}$-boundary passes between the $B_{g^{-k+1}x_0}$- and 
the $B_{g^{-k}x_0}$-boundaries and then turns around the $B_{g^{-j}x_0}$-boundary
for $j<k$ and then connects to the $A_{g^{-k}x_0}$-boundary parallel to the other $A_x$-edges.

We now want to transform this marking into the standard marking by a finite sequence of LTG-moves. To do so, recall
that for every $k>l$ the marking is of the form

\be
\raisebox{20pt}{
\bp(120,65)
\put(0,0){\pic{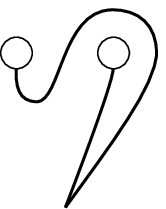}}
\put(2,52){\scriptsize $A^{-k}$}
\put(29,60){\scriptsize $B^{-l}$}
\ep
}
\ee
Hence we need to apply the LTG-move $\B_{B_{g^{-l}x_0},A_{g^{-k}x_0}}$ to turn this into the marking

\be
\raisebox{20pt}{
\bp(100,65)
\put(0,0){\pic{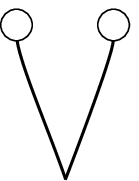}}
\put(2,53){\scriptsize $A^{-k}$}
\put(29,53){\scriptsize $B^{-l}$}
\ep
}
\ee
Now starting with the $B_{g^{-n+2}x_0}$-boundary, for $l=n-2,\dots,1$ we successively apply $\B_{B_{g^{-l}x_0},A_{g^{-k}x_0}}$ for 
$k=l+1,\dots,n-1$. 
This finally gives the standard marking on $S_{2n+2}$. 
To translate these LTG-moves into a morphism in $\C$, we introduce auxiliary morphisms $f^{(k)}$ for $k=n,\dots, 1$ with
\be
\raisebox{20pt}{
\bp(21,50)
\put(0,-10){\rib{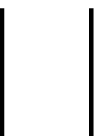}}
\put(-40,5){$f^{(n)}=$}
\put(-8,31){\scriptsize $A^{-n+1}$}
\put(-8,-19){\scriptsize $A^{-n+1}$}
\put(21,31){\scriptsize $B^{-n+1}$}
\put(21,-19){\scriptsize $B^{-n+1}$}
\ep
}
\ee
the identity and $f^{(k)}$ for $k=n-1,\dots, 1$ recursively
\be
\raisebox{20pt}{
\bp(110,110)
\put(0,0){\rib{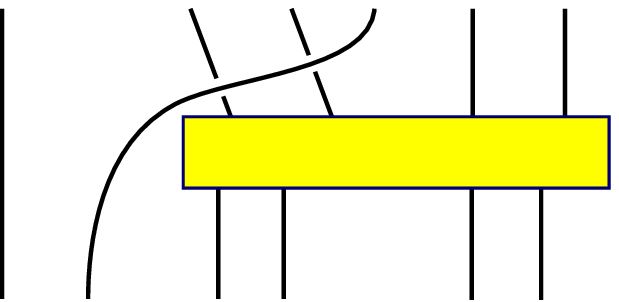}}
\put(-50,40){$f^{(k)}=$}
\put(105,40){\scriptsize $f^{(k+1)}$}
\put(-7,-10){\scriptsize $A^{-k+1}B^{-k+1}$}
\put(55,-10){\scriptsize $A^{-k}$}
\put(78,-10){\scriptsize $B^{-k}$}
\put(108,-10){\scriptsize $\dots$}
\put(126,-10){\scriptsize $A^{-n+1}B^{-n+1}$}
\put(-7,88){\scriptsize $A^{-k+1}$}
\put(47,88){\scriptsize $A^{-k}$}
\put(60,88){\scriptsize $\dots$}
\put(74,88){\scriptsize $A^{-n+1}$}
\put(100,88){\scriptsize $B^{-k+1}$}
\put(131,88){\scriptsize $B^{-k}$}
\put(144,88){\scriptsize $\dots$}
\put(158,88){\scriptsize $B^{-n+1}$}
\ep
}
\ee
So the step $f^{(l+1)}\to f^{(l)}$ resembles the application of the LTG-moves $\B_{B_{g^{-l+1}x_0},A_{g^{-k}x_0}}$ for 
$k=l,\dots,n-1$. 

Alltogether this gives the associativity constraint
\be
\psi_{(A_x),(B_x),M}=f^{(1)}\otimes\id_M
\ee

An example of this morphism in the case $n=4$ is depicted in 

\be\label{eqn:psiexample}
\raisebox{20pt}{
\bp(110,110)
\put(0,0){\rib{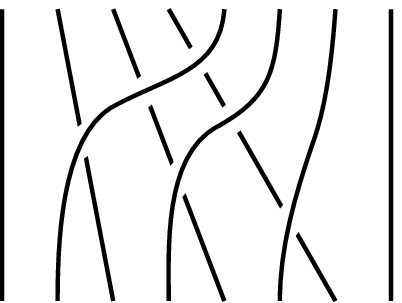}}
\put(-90,40){$\psi_{(A_x),(B_x),M}=$}
\put(-4,-10){\scriptsize $A^{0}$}
\put(12,-10){\scriptsize $B^{0}$}
\put(28,-10){\scriptsize $A^{-1}$}
\put(44,-10){\scriptsize $B^{-1}$}
\put(61,-10){\scriptsize $A^{-2}$}
\put(78,-10){\scriptsize $B^{-2}$}
\put(93,-10){\scriptsize $A^{-3}$}
\put(109,-10){\scriptsize $B^{-3}$}
\put(-4,90){\scriptsize $A^{0}$}
\put(12,90){\scriptsize $A^{-1}$}
\put(28,90){\scriptsize $A^{-2}$}
\put(44,90){\scriptsize $A^{-3}$}
\put(61,90){\scriptsize $B^{-0}$}
\put(78,90){\scriptsize $B^{-1}$}
\put(93,90){\scriptsize $B^{-2}$}
\put(109,90){\scriptsize $B^{-3}$}
\ep
}
\ee

For general $n$, the morphism $\psi_{(A_x),(B_x),M}$ is a shuffle that moves all $B$-objects over the $A$-objects.

By the general arguments of \cite{kpGMF}, these associativity constraints satisfy the mixed pentagon axiom, which can also easily be verified
by hand. We summarize our findings for arbitrary groups in the following
\begin{thm}\label{thm:firstmain2}
For any finite group $G$, any finite $G$-set $\X$ and any $g\in G$ the functor
\be
\bearll
\C^{\boxtimes\X}\times \C^{\boxtimes O_g}&\rightarrow \C^{\boxtimes O_g}\\
(A_x)_{x\in\X}\times (M_o)_{o\in O_g}&\mapsto (A_{x_o}A_{g^{-1}x_o}\cdots A_{g^{-|o|+1}x_o}M_o)_{o\in O_g}
\eear
\ee
with $x_o$ the smallest element in the $\la g \ra$-orbit $o$, together with the associativity constraints
\be
\Psi_{A,B,M}=(\psi^o_{(A_x)_{x\in o},(B_x)_{x\in o},M})_{o\in O_g}\,,
\ee
with the morphisms $\psi^o$ as in equation \erf{eqn:psiexample}, endows $\C^{\boxtimes O_g}$ with the structure of a module category over the tensor category $\C^{\boxtimes \X}$.
\end{thm}

\section{Permutation modular invariants}\label{sec:invariants}

For any semisimple module category $(\mathcal{M},\otimes,\psi)$ over a braided fusion 
category $\mathcal{D}$ the category $\Endfun_{\mathcal{D}}(\mathcal{M})$ of $\D$-module endofunctors of $\M$ is again a monoidal category that acts on $\M$.
Recall from \cite[Section 5.1,5.2]{ostrik} the following definition:
\begin{Def}\label{def:alphaind}
Define two functors
\be
\alpha^\pm:\mathcal{D}\to \Endfun_{\mathcal{D}}(\mathcal{M})
\ee
by putting 
\be
\alpha^\pm(U)(M):=U\otimes M
\ee
as functors and with the following module functor constraints for $\alpha^\pm(U)$:
\be\label{eqn:alphaconstraints}
\bearl
\gamma^{U,+}_{V,M} := \psi^{}_{V,U,M} \circ (c_{U,V}^{}\oti\id_M) \circ \psi^{-1}_{U,V,M} 
\qquad{\rm and} \\{}\\[-.8em]
\gamma^{U,-}_{V,M} := \psi^{}_{V,U,M} \circ (c_{U,V}^{-1}\oti\id_M) \circ \psi^{-1}_{U,V,M} \,,
\eear
\ee
where $c_{U,V}$ is the braiding in $\D$. The functors $\alpha^\pm$ are called the \emph{$\alpha$-induction functors}.
\end{Def}

\begin{rem}
The associativity constraints $\psi$ of the module category $\mathcal{M}$ endow the functors $\alpha^\pm$ with the structure of monoidal functors.
\end{rem}

If $U$ is an object of $\D$, we abbreviate $\alpha^\pm(U)\equiv\alpha^\pm_U$; if $U_k$ is a simple object of $\mathcal{D}$, we write $\alpha^\pm(U_k)\equiv\alpha^\pm_k$.
Now for any two simple objects $U_i$, $U_j$ of $\mathcal{D}$ we define the non-negative integers
\be
Z_{i,j}:=\dimk\hom_{\Endfun_{\mathcal{D}}(\mathcal{M})}(\alpha^+_i,\alpha^-_j)
\ee
The $|\I|\times|\I|$-matrix $Z(\mathcal{M}/\mathcal{D}):=(Z_{i,j})$ then obeys the requirements on a modular invariant (see \cite[Theorem 5.1]{TFT1}).
It is the aim of this section to understand the structure of this matrix in the case that the module category
under consideration is given by the data of the $G$-equivariant modular functor $\tauzz$.

So we fix an element $g\in G$ and examine the module category $\C^\X_g=\C^{\boxtimes O_g}$ over $\C^\X_1=\C^{\boxtimes\X}$. 
We will continue in two steps: First we show that certain entries of the matrix $Z(\CX_g/\CX_1)$ are non-zero. Then
we show that $Z(\CX_g/\CX_1)$ is a permutation matrix, i.e. it contains precicely one entry $1$ in every row and column
and $0$ elsewhere. This already fixes the whole matrix.

At first we will give for every object $U$ in $\C^{\boxtimes\X}$ an invertible natural transformation
\be
\Gamma^U:\alpha^+(U)\Rightarrow\alpha^-({}^gU)
\ee
between module functors.
For $U=U_{\bar\imath}$ a simple object, this will be a non-zero element in 
\be
\hom_{\Endfun_{\C^{\boxtimes\X}}(\C^{\boxtimes O_g})}(\alpha^+_{\bar\imath},\alpha^-_{g\bar\imath})
\ee
so that this space is non-zero. As in section \ref{sec:main} we can restrict our discussion to only one factor of
$\C^{\boxtimes O_g}$ by giving $\Gamma^U$ for every factor separately
so that we are again in the situation of convention \ref{conv:cyclic}. Hence we will consider the problem where
$\C^{\boxtimes\X}$ acts on a single copy of $\C$.

From now on we fix the object $U=(U_x)_{x\in \X}$ in
$\C^{\boxtimes\X}$ and write $\Gamma$ instead of $\Gamma^U$. For $M$ in $\C$ we find  
\be
\bearll
\alpha^+_{U}(M)&=U_{x_o}U_{g^{-1}x_o}\cdots U_{g^{-n+1}x_o}M\\
\alpha^-_{{}^gU}(M)&=U_{g^{-1}x_o}\cdots U_{g^{-n+1}x_o}U_{x_o}M
\eear
\ee

Hence we let
\be
\raisebox{20pt}{
\bp(70,110)
\put(0,0){\rib{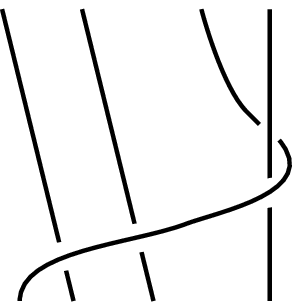}}
\put(-50,40){$\Gamma_M:=$}
\put(0,-10){\scriptsize $U^{0}$}
\put(18,-10){\scriptsize $U^{-1}$}
\put(33,-10){\scriptsize $\dots$}
\put(45,-10){\scriptsize $U^{-n+1}$}
\put(75,-10){\scriptsize $M$}
\put(-3,90){\scriptsize $U^{-1}$}
\put(22,90){\scriptsize $U^{-n+1}$}
\put(9,90){\scriptsize $\dots$}
\put(55,90){\scriptsize $U^{0}$}
\put(75,90){\scriptsize $M$}
\ep
}
\ee
In formulas:
\be
\Gamma_M=\left[\id_{U_{g^{-1}x_o}\cdots U_{g^{-n+1}x_o}}\otimes (c_{M,U_{x_o}}\circ c_{U_{x_o},M})\right]\circ \left[c_{U_{x_o},U_{g^{-1}x_o}\cdots U_{g^{-n+1}x_o}}\otimes\id_M\right]
\ee

\begin{lemma}\label{lem:nattrans}
$\Gamma$ is a non-zero natural transformation between the module functors $\alpha^+_{U}$ and $\alpha^-_{{}^gU}$.
\end{lemma}
\begin{proof}
Obviously $\Gamma_M$ is natural in $M$ and invertible, hence non-zero. To show that $\Gamma$ is a natural transformation
of module functors, we have to show that for another object $V=(V_x)$ of $\C^{\boxtimes \X}$ the following
compatibility with the module functor constraints \erf{eqn:alphaconstraints} holds:
\be
(\id_{(V_x)}\otimes \Gamma_M)\circ \gamma^{U,+}_{V,M}=\gamma^{{{}^gU},-}_{V,M}\circ \Gamma_{(V_x)\otimes M}
\ee
Spelling out all occurring morphisms, this amounts to
\be\label{eq:ind}
\bearl
\left[\id_{V_{x_o}\cdots V_{g^{-n+1}x_o}}\otimes\Gamma_M\right]\circ\psi_{(V_x),(U_x),M}\circ\left[c^\X_{(U_x),(V_x)}\otimes \id_M\right]\circ\psi_{(U_{x}),(V_x),M}^{-1}=\\
\psi_{(V_x),{}^g(U_{x}),M}\circ\left[(c^\X_{(V_x),{}^g(U_{x})})^{-1}\otimes\id_M\right]\circ\psi_{{}^g(U_{x}),(V_x),M}^{-1}\circ\Gamma_{V_{x_o}\cdots V_{g^{-n+1}x_o}M}
\eear
\ee
where $c^\X_{(U_x),(V_x)}$ denotes the braiding of two objects in $\C^{\boxtimes\X}$. This is an equality in
$$
\homc(U_{x_o}\cdots U_{g^{-n+1}x_o}V_{x_o}\cdots V_{g^{-n+1}x_o}M,V_{x_o}\cdots V_{g^{-n+1}x_o}U_{g^{-1}x_o}\cdots U_{g^{-n+1}x_o}U_{x_o}M)\,\,.
$$
Denote the left hand side 
of \erf{eq:ind} by $L^n(U_{x_o},\dots,U_{g^{-n+1}x_o};V_{x_o},\dots,V_{g^{-n+1}x_o};M)$ and the right hand side by 
$R^n(U_{x_o},\dots,U_{g^{-n+1}x_o};V_{x_o},\dots,V_{g^{-n+1}x_o};M)$. The endomorphism
\be
\raisebox{20pt}{
\bp(230,110)
\put(0,0){\rib{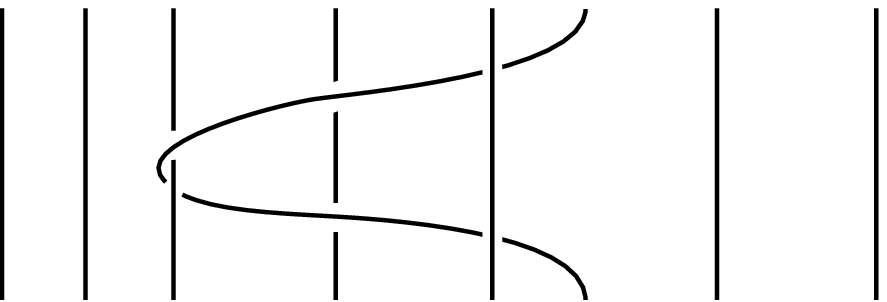}}
\put(-50,40){$F=$}
\put(-3,-10){\scriptsize $U^{0}$}
\put(22,-10){\scriptsize $U^{-1}$}
\put(22,90){\scriptsize $U^{-1}$}
\put(-3,90){\scriptsize $U^{0}$}
\put(68,-10){\scriptsize $(U^{-3}\dots U^{-n+1})$}
\put(68,90){\scriptsize $(U^{-3}\dots U^{-n+1})$}
\put(47,-10){\scriptsize $U^{-2}$}
\put(47,90){\scriptsize $U^{-2}$}
\put(139,-10){\scriptsize $V^{0}$}
\put(139,90){\scriptsize $V^{0}$}
\put(163,-10){\scriptsize $V^{-1}$}
\put(163,90){\scriptsize $V^{-1}$}
\put(182,-10){\scriptsize $(V^{-2}\dots V^{-n+1})$}
\put(182,90){\scriptsize $(V^{-2}\dots V^{-n+1})$}
\put(250,-10){\scriptsize $M$}
\put(250,90){\scriptsize $M$}
\ep
}
\ee
is obviously invertible and an easy but lengthy graphical calculation shows that it obeys
\be
\bearl
L^n(U_{x_o},\dots,U_{g^{-n+1}x_o};V_{x_o},\dots,V_{g^{-n+1}x_o};M)\circ F=\\
L^{n-1}(U_{x_o},U_{g^{-1}x_o}U_{g^{-2}x_o},\dots,U_{g^{-n+1}x_o};V_{x_o},V_{g^{-1}x_o}V_{g^{-2}x_o},\dots,V_{g^{-n+1}x_o};M)
\eear
\ee
and
\be
\bearl
R^n(U_{x_o},\dots,U_{g^{-n+1}x_o};V_{x_o},\dots,V_{g^{-n+1}x_o};M)\circ F=\\
R^{n-1}(U_{x_o},U_{g^{-1}x_o}U_{g^{-2}x_o},\dots,U_{g^{-n+1}x_o};V_{x_o},V_{g^{-1}x_o}V_{g^{-2}x_o},\dots,V_{g^{-n+1}x_o};M)
\eear
\ee
Hence \erf{eq:ind} holds by induction,
the case $n=2$ is an easy calculation using only relations in the braid group on five strands.
\end{proof}

When we apply lemma \ref{lem:nattrans} to a simple object $U_{\bar\imath}=(U_i)_{i\in\bar\imath}$ with 
$\bar\imath\in \I^\X$ we see
\be
Z_{\bar\imath,g\bar\imath}=\dimk\hom_{\Endfun_{\C^{\boxtimes\X}}(\C)}(\alpha^+_{\bar\imath},\alpha^-_{g\bar\imath})\neq 0\,.
\ee

\begin{Def}
An {\em Azumaya category} $\M$ over a braided monoidal category $\D$
is a left module category $\M$ over $\D$ for which the two monoidal functors
$\alpha^\pm$ from $(\D,\otimes)$ to $(\Endfun_\D(\M),\circ)$ are equivalences.
\\[1pt]
An {\em Azumaya algebra\/} $A$ in a braided monoidal category $\D$
is an algebra $A$ in $\D$ such that the category of right $A$-modules is an
Azumaya category over $\D$.
\end{Def}
\begin{rem}
This definition is equivalent to the definition given in \cite[Section 3]{vazh}.
\end{rem}

By \cite[Theorem 6.1]{ENO} every bimodule category $\M$ over $\D$, which is part of an equivariant monoidal category
with neutral component $\D$, is invertible with respect to a tensor product of module categories. We will need the following
criterion for invertibility of a module category. We assume that $\M$ is a module category that is turned into a
bimodule category by using the braiding of $\D$.

\begin{lemma}\label{lem:invAzu}
A semisimple module category $\mathcal{M}$ over a modular category $\mathcal{D}$ is invertible if and only if
it is equivalent to $A\!-\!{\rm mod}$, as a module category over $\mathcal{D}$, for some Azumaya algebra $A$ in $\mathcal{D}$.
\end{lemma}
\begin{proof}
By proposition 4.2 and section 5.4 of \cite{ENO} invertibility of $\M$ is equivalent to $\M$ being an Azumaya category.
If $\M$ is invertible, by \cite[Corollary 4.4]{ENO} it is indecomposable over $\D$. It follows
from \cite[Theorem 1]{ostrik} that as a module category $\M$ is equivalent to $A\!-\!{\rm mod}$ for some algebra $A$ in $\D$,
which then is Azumaya.
\end{proof}

\void{

The proof of \cite[Theorem 1]{ostrik} uses the formalism of internal homs. As we will also use this in the following 
discussion we recall the definition from \cite[Section 3.2]{ostrik}.

\begin{Def}
Let $\M$ be a semisimple (left) module category over a semisimple fusion category $\D$. For any two objects $M_1$ and $M_2$ of $\M$
 their internal Hom $\underline\hom (M_1,M_2)$ is
an object of $\D$ that represents the functor
\be
\begin{array}{ll}
\D&\to\vectk\\
U&\mapsto \hom_\M(U\otimes M_1,M_2).
\end {array}
\ee

\end{Def}

The internal Hom always exists, when $\D$ and $\M$ are semisimple and $\D$ has only finitely many isomorphism classes of simple objects.
For $M=M_1=M_2$ we write $\intEnd(M):=\underline\hom(M,M)$ and call this object the internal End of $M$.
There is a canonical evaluation morphism $e_{M_1,M_2}:\underline\hom(M_1,M_2)\otimes M_1\to M_2$ in $\M$ that induces
an associative composition morphism 
\be
\underline\hom(M_2,M_3)\underline\hom(M_1,M_2)\to\underline\hom(M_1,M_3).
\ee
This morphism gives the structure of an associative algebra on the internal End of an object $M$ in $\M$.

}


\begin{thm}
The modular matrix $Z(\CX_g/\CX_1)$ for the module category described in theorem \ref{thm:firstmain2} reads
\be
Z(\CX_g/\CX_1)_{\bar \imath,\bar \jmath}=\delta_{\bar \jmath,g\bar\imath}
\ee
where $\bar\imath, \bar\jmath\in\I^\X$ label the simple objects of $\CX_1=\C^{\boxtimes\X}$ and 
$g\bar\imath$ is the multi-index $\bar\imath$ permuted by the action of the group element $g\in G$.
\end{thm}
\begin{proof}
Since $\tauzz$ is a $G$-equivariant modular functor, it induces (\cite{kpGMF}) the structure of a $G$-equivariant category on
$\bigoplus \CX_h$. The module category $\CX_g$ over $\CX_1$ is part of this larger structure, by \cite[Theorem 6.1]{ENO} 
it is invertible, hence by lemma \ref{lem:invAzu} equivalent to the category $A\!-\!{\rm mod}$
for some Azumaya algebra $A$ in $\CX_1$. As $A$ is Azumaya, the functors $\alpha^\pm$ are equivalences
$$
\alpha^\pm:\,\CX_1\to \Endfun_{\CX_1}(\CX_g)\cong A\!-\!{\rm bimod}
$$
of tensor categories. Hence $A\!-\!{\rm bimod}$ is semisimple and for a simple object $U_{\bar\imath}$ in $\CX_1$ the objects 
$\alpha^\pm (U_{\bar\imath})$ in $A\!-\!{\rm bimod}$ are again simple. By semisimplicity of $A\!-\!{\rm bimod}$,
the matrix 
$$
Z(\CX_g/\CX_1)_{\bar\imath,\bar\jmath}:=\dimk\hom_{\Endfun_{\CX_1}(\CX_g)}(\alpha^+_{\bar\imath},\alpha^-_{\bar\jmath})
$$
has exactly one entry $1$ in every row and every column and $0$ elsewhere. 
By lemma \ref{lem:nattrans} we find that in every row and column the numbers $Z_{\bar\imath,g\bar\imath}$ are non-zero. 
\end{proof}

\newpage
\bibliography{bib}{}
\bibliographystyle{alpha}
\end{document}